\documentclass[english]{article}
\usepackage[T1]{fontenc}
\usepackage[latin1]{inputenc}
\usepackage{esint}

\makeatletter

\newtheorem{theorem}{Theorem}

\newtheorem{lemma}[theorem]{Lemma}

\newtheorem{remark}[theorem]{Remark}

\date{}

\makeatother

\makeatother

\usepackage{babel}\makeatother

\makeatother

\usepackage{babel}

\begin{document}

\title{Accumulation on the boundary for one-dimensional stochastic particle
system}

\author{V. A. Malyshev, A. A. Zamyatin}
\maketitle
\begin{abstract}
We consider infinite particle system on the positive half-line moving
independently of each other. When a particle hits the boundary it
immediately disappears, and the boundary moves to the right on some
fixed quantity (particle size). We study the speed of the boundary
movement (growth). Possible applications - dynamics of the traffic
jam growth, growth of thrombus, epitaxy. Nontrivial mathematics is
related to the correlation between particle dynamics and boundary
growth. 
\end{abstract}

\section{Introduction}

Hitting the boundary by the random walk is a classical problem of
probability theory. More complicated is the problem of hitting the
moving boundary, see \cite{Shir}. Here we consider the case when
the movement of the boundary is correlated with the movement of the
particles. Namely, we consider infinite number of particles, performing
random walks on the positive half-line, which stick to the boundary
when hit it. Moreover, when a particle hits the boundary, the boundary
moves to the right on some quantity, {}``size'' of the particle.
That is the boundary moves due to accumulation of particles on it.
As far as we know, such problems were not considered earlier.

Similar problems one can often encounter in applications - growth
of thrombus, epitaxy and other methods of coating the surface with
metal particles. In many cases the effect of correlation between the
movements of particles and the boundary can be neglected. Here, on
the contrary, we consider the influence of this correlation on the
growth speed of the boundary. Such effects can be observed while quick
formation of traffic jams.

In the paper we present two results. First one for nonzero drift of
the particles, second - for zero drift. In the first case one can
find exact asymptotic rate of growth, in the second case one can get
only the order of growth.

\section{Formulation of the problem and the results}

We consider infinite system of particles on the half-line, which,
while hitting the boundary, adhere to it (accumulate on it). As a
result the boundary moves to the right. Now we give exact formulation.

Assume that at time $0$ on the half-line $R_{+}$the random configuration
of particles at the points \[
0\leq x_{1}=x_{1}(0)<...<x_{n}=x_{n}(0)<...\]
 is distributed as the point Poisson process with density $\lambda$.

Each particle, before hitting the boundary, moves independently of
other particles as\begin{equation}
x_{i}(t)=x_{j}-vt+w_{j}(t)\label{eq:mv}\end{equation}
 where $v\geq0$ is the constant drift, and $w_{i}(t)$ is the standard
Wiener process with zero mean. At the hitting moment the particles
are absorbed by the boundary and disappear, moreover the boundary
moves to the right on some $\delta$, the {}``particle size''.

More exactly, the boundary which at time $t=0$ coincides with point
$0\in R_{+}$, also moves. The law $\xi(t)$ of its movement models
the sticking of particles to the boundary. Namely, $\xi(0)=0$, and
$\xi(t)$ is a non-decreasing point-wise constant function with jumps\[
\xi(t_{j}+0)=\xi(t_{j})+\delta k_{j}\]
 at random time moments\[
0\leq t_{1}<...<t_{j}<...\]
 where $\delta>0$ is some constant, the {}``particle size'', the
moments $t_{j}$ and the non-negative integers $k_{j}$ being defined
by the following recurrent procedure. Let $t_{1}$ be the first instant,
when one of the particles hits $0$, then $k_{1}$ is the minimal
integer $k>0$ such that in the interval $(0,k\delta]$ at this moment
$t_{1}$ there are not more than $k-1$ particles. For example, $k=1$
if and only if at time $t_{1}$ there are no particles inside the
interval $(0,\delta]$. Further on by induction, $t_{n+1}$ is defined
as the first moment when one of the remaining particles hits the point
$(k_{1}+...+k_{n})\delta$, and $k_{n+1}$ is defined correspondingly
as the minimal integer $k>0$ such that inside the interval $(\xi(t_{n}),k\delta]$
at time $t_{n+1}$ there is no more than $k-1$ particles.

Below the following statement will be proved.

\begin{lemma}\label{lemma0}

If $\lambda<\delta^{-1}$, then with probability 1 all $k_{i}$ and
$t_{i}$ are finite, and moreover $t_{i}\to\infty$. On the contrary,
if $\lambda>\delta^{-1}$, the boundary reaches infinity for finite
time with positive probability.

\end{lemma}

Now we can formulate the main result. Let $N=N(t)$ be the number
of particles, absorbed by the boundary to time $t$, then $\xi=\xi(t)=\delta N(t)$
is the coordinate of the particle at time $t$.

\begin{theorem}

Let $\lambda<\delta^{-1}$. If $v>0$, the particle movement is asymptotically
linear, that is as $t\to\infty$ a. s. \[
\frac{\xi(t)}{t}\to V=v\frac{\delta\lambda}{1-\delta\lambda}\]
 If $v=0$, then for $t$ sufficiently large \[
C_{0}\leq\frac{MN(t)}{\sqrt{t}}\leq C\]
 for some constants\[
\lambda\sqrt{\frac{2}{\pi}}<C_{0}<C<\infty\]

\end{theorem}

\begin{remark}We will show later that under the condition that the
boundary does not move at all and $v=0,$ then the mean number of
particles hitting the boundary during time $t$ equals $\lambda\sqrt{\frac{2t}{\pi}}.$
From the second assertion of the theorem it follows in particular
that the mean number of particles hitting the boundary is asymptotically
greater. We conjecture that $\frac{MN(t)}{\sqrt{t}}$ tends to some
constant as $t\to\infty$, and moreover the distribution of the random
variable $\frac{N(t)}{\sqrt{t}}$ tends to some continuous distribution
as $t\to\infty$.

\end{remark}

\section{Proofs}

\paragraph{Proof of Lemma \ref{lemma0}}

For a given configuration $X=\{x_{i}(0)\}$ denote $p(X,n)$ the probability
that the first particle hitting the boundary is the particle $x_{n}(0)$.

To prove the first part of the Lemma introduce the following auxiliary
model. At time $t=0$ for a given configuration on the positive axis
we add Poisson configuration on the negative axis $Y$, having the
same density $\lambda.$ Then we will get Poisson configuration on
all line $R$, all these particles have red color. Moreover, there
is one more (blue) particle, which we put with probability $p(X,n)$
at the point $x_{n}(0)>0.$ The particles move according to (\ref{eq:mv}),
independently of each other. At time $t=0$ the boundary is at the
point $\xi'(0)=0$. The random variable $t'_{1}$ is defined as the
first moment when the blue particle hits $0$, and $k'_{1}$ are defined
as earlier. The difference from the basic model is that the boundary
does not influence other particles at all. Thus at any moment $t$
the red particles $y_{i}(t),-\infty<j<\infty,$ have Poisson field
distribution (Doob-Dobrushin theorem, \cite{Dobr}). Moreover at random
time $t'_{1}$ the red particles have Poisson field distribution.

Let $\xi'(t)$ be the position of the boundary at time $t.$ Then
for all $t$ a. s. $\xi'(t)\geq\xi(t).$ In particular, $k'_{1}\geq k_{1}$
a. s. In fact, let $\Omega$ be the probability space of the auxiliary
model. The basic model is defined on the same probability space, as
it can be obtained from the auxiliary model by deleting the glue particle
and the particles situated to the left of $0$. Hence in the auxiliary
model the number of collisions with the boundary is not less than
in the basic model.

Let us prove that $k'_{1}$ is finite a. s. for any $n$. The random
variable $k'_{1}$ is defined by the stopping time of the following
Markov chain $m_{1},...,m_{i},...$ with discrete time and state space
$Z_{+}$. The state $0$ is an absorbing state. Denote $m_{1}$the
number of particles in the interval $[0,\delta]$ at time $t_{1}$,
and $m_{2}$ the number of particles in the interval $[\delta,(m_{1}+1)\delta]$,
if $m_{1}>0$, $m_{3}$ the number of particles in the interval $[(m_{1}+1)\delta,(m_{1}+1)\delta+m_{2}\delta],$
if $m_{1}>0$, $m_{2}>0$, and so on. The sequence $m_{i}$ constitute
time homogeneous Markov chain with transition probabilities\[
p_{m_{i}m_{i+1}}=\frac{(m_{i}\lambda\delta)^{m_{i+1}}}{m_{i+1}!}e^{-m_{i}\lambda\delta}\]
 The boundary stops when $m_{i}=0.$ One step increment is \[
M(m_{i+1}-m_{i}|m_{i}=k)=k(\lambda\delta-1)\]
 The chain hits $0$ with probability $1$, if $\lambda\delta-1<0$.

By induction, using similar trick with appending a new blue particle
at any moments $t'_{i}$, the finiteness of all $k'_{i}$ follows.

The second statement of the Lemma is not used in the proof of the
Theorem, and it will be convenient to prove it at the end of the paper
after the proof of the Theorem.

\paragraph{Idea of the proof of the Theorem}

One can write down two equalities relating random variables $\xi(t)$
and $N(t)$: \begin{equation}
\xi(t)=N(t)\delta\label{equa1}\end{equation}
 \begin{equation}
N(t)=n_{0}(t,\xi(t))+n_{1}(t,\xi(t))\label{equa2}\end{equation}
 where $n_{0}(t,\xi)$ is the number of particles which were at time
$0$ inside the interval $[0,\xi]$ and were absorbed by the boundary
before time $t$, $n_{1}(t,\xi)$ is the number of particles which
were at time $0$ to the right of the point $\xi$ and were absorbed
by the boundary before time $t$. We can write the following equation
\[
\delta^{-1}\xi(t)=n_{0}(t,\xi(t))+n_{1}(t,\xi(t))\]
 It will appear useful if we shall have good expressions for $n_{0}$
and $n_{1}$.

The first result of the theorem can be explained for the case when
for nonzero drift there is no fluctuations. More exactly, the particles,
situated initially at the points $\lambda^{-1}k$ with integer $k$,
move with constant velocity $v$. Then from equations (\ref{equa1}),(\ref{equa2})
we have\begin{equation}
\xi=N\delta,N=[\xi\lambda]+[\lambda vt]\end{equation}
 Hence\[
N-[\lambda N\delta]=[\lambda vt]\]
 \[
N=N(t)=\frac{[\delta\lambda vt]}{1-\frac{1}{N}[N\delta\lambda]}\sim\frac{\delta\lambda vt}{1-\delta\lambda}\]
 for large $t$. Then\[
\frac{\xi(t)}{t}\to V=v\frac{\delta\lambda}{1-\delta\lambda}\]

\subsection{The case of nonzero drift}

Fix arbitrary $\epsilon>0,$ so that $(\lambda+\epsilon)\delta<1.$
Let $c_{\epsilon}=(1-(\lambda+\epsilon)\delta)^{-1}(\lambda+\epsilon)(v+\epsilon).$
The event $A(t)=\{N(t,\omega)>c_{\epsilon}t\}$ can be written as
the union

\[
A(t)=\cup_{n}B(n,t)\]
 where $B(n,t)$ is the event that to the time $t$ the boundary has
absorbed exactly $n$ particles, situated initially at the interval
$(0,c_{\epsilon}t\delta+(v+\epsilon)t)$ and more than $m=c_{\epsilon}t-n$
particles, situated initially at the interval $(c_{\epsilon}t\delta+(v+\epsilon)t,\infty)$.
Denote $B_{\epsilon}(t)$ the union of events $B(n,t)$ with $n\geq(\lambda+\epsilon)(\delta c_{\epsilon}t+(v+\epsilon)t)=c_{\epsilon}t$.
The event $B_{\epsilon}(t)$ has exponentially small probability,
as it belongs to the event $C_{\epsilon}(t)$ that at time $0$ in
the interval $(0,c_{\epsilon}t\delta+(v+\epsilon)t)$ there were not
less than $(\lambda+\epsilon)(\delta c_{\epsilon}t+(v+\epsilon)t)=c_{\epsilon}t$
particles. The probability of the latter event can be estimated via
Poisson distribution:\begin{equation}
P(C_{\epsilon}(t))\leq C_{0}e^{-h(\epsilon)t}\label{eq:est1}\end{equation}
 where constant $C_{0}$ depends on $\lambda$ and $\epsilon,$ $h(\epsilon)>0.$

Denote $D_{\epsilon}(t)$ the union of the events $B(n,t)$ with $n<(\lambda+\epsilon)(\delta c_{\epsilon}t+(v+\epsilon)t)=c_{\epsilon}t$.
The event $D_{\epsilon}(t)$ belongs to the event that at least one
of the particles, situated initially in $(c_{\epsilon}t\delta+(v+\epsilon)t,\infty)$,
before $t$ will be absorbed by the boundary. The latter event in
its turn belongs to the event $J_{\epsilon}(t),$ that at least one
of the particles, situated initially in $(c_{\epsilon}t\delta+(v+\epsilon)t,\infty),$
had been to the left of the point $c_{\epsilon}t\delta$ at least
once before $t$. Thus,\[
A(t)\subset C_{\epsilon}(t)\cup J_{\epsilon}(t)\]

The probability of the event $J_{\epsilon}(t)$ can be estimated from
above by the expectation of the number of particles, situated initially
in $(c_{\epsilon}t\delta+(v+\epsilon)t,\infty),$ and visiting at
least once the left side of the point $c_{\epsilon}t\delta$ before
$t$ \[
P(J_{\epsilon}(t))\leq\lambda\int_{c_{\epsilon}t\delta+(v+\epsilon)t}^{\infty}p(t,x)dx\]
 where $p(t,x)$ is the probability that the particle started at $x$,
$x>c_{\epsilon}t\delta+(v+\epsilon)t,$ during $t$ will visit the
left side of the point $c_{\epsilon}t\delta.$ More exactly, \begin{equation}
p(t,x)=P(\min_{s\leq t}(x-vs+w(s))\leq c_{\epsilon}t\delta)\label{eq:p1}\end{equation}
 From the definition of Wiener process it follows that \begin{eqnarray*}
P(\min_{s\leq t}(x-vs+w(s))\leq c_{\epsilon}t\delta) & = & P(\min_{s\leq t}(-vs+w(s))\leq-x+c_{\epsilon}t\delta)\\
 & = & P(\max_{s\leq t}(vs+w(s))\geq x-c_{\epsilon}t\delta)\end{eqnarray*}
 To find $p(t,x)$ we will use the well known formula for the distribution
of the maximum of Wiener process with drift (see for example, \cite{Shir},
p. 935) \begin{equation}
P(\max_{s\leq t}(vs+w(s))\geq x)=1-\Phi\left(\frac{x-vt}{\sqrt{t}}\right)+e^{2vx}\Phi\left(\frac{-x-vt}{\sqrt{t}}\right)\label{eq:distr_max}\end{equation}
 for $\Phi(z)=\frac{1}{\sqrt{2\pi}}\int_{-\infty}^{z}e^{-y^{2}/2}dy.$
Then\begin{eqnarray}
P(J_{\epsilon}(t)) & \leq & \lambda\int_{c_{\epsilon}t\delta+(v+\epsilon)t}^{\infty}\bar{\Phi}\left(\frac{x-c_{\epsilon}t\delta-vt}{\sqrt{t}}\right)dx\nonumber \\
 &  & +\lambda\int_{c_{\epsilon}t\delta+(v+\epsilon)t}^{\infty}e^{2vx}\Phi\left(\frac{-x+c_{\epsilon}t\delta-vt}{\sqrt{t}}\right)dx\label{eq:sum}\\
 & = & I_{1}(\epsilon,t)+I_{2}(\epsilon,t)\nonumber \end{eqnarray}
 where $\bar{\Phi}=1-\Phi.$ Change of variable $z=\frac{x-c_{\epsilon}t\delta-vt}{\sqrt{t}}$
in the first integral in the formula (\ref{eq:sum}) gives\[
I_{1}(\epsilon,t)=\lambda\sqrt{t}\int_{\epsilon\sqrt{t}}^{\infty}\bar{\Phi}(z)dz\leq\lambda\sqrt{\frac{t}{2\pi}}e^{-\frac{\epsilon^{2}t}{2}}\]
 with $\epsilon\sqrt{t}>1.$

Change of variables $z=\frac{x-c_{\epsilon}t\delta+vt}{\sqrt{t}}$
in the second integral in (\ref{eq:sum}) gives\begin{eqnarray*}
I_{2}(\epsilon,t) & = & \lambda\sqrt{t}e^{-2v^{2}t}\int_{(2v+\epsilon)\sqrt{t}}^{\infty}e^{2v\sqrt{t}z}\Phi(-z)dz\\
 & = & \lambda\sqrt{t}e^{-2v^{2}t}\int_{(2v+\epsilon)\sqrt{t}}^{\infty}e^{2v\sqrt{t}z}\bar{\Phi}(z)dz\end{eqnarray*}
 For $\epsilon\sqrt{t}>1$ we have\begin{eqnarray*}
I_{2}(\epsilon,t) & \leq & \lambda\sqrt{\frac{t}{2\pi}}\int_{(2v+\epsilon)\sqrt{t}}^{\infty}e^{-(z^{2}/2-2vz\sqrt{t}+2v^{2}t)}dz\\
 & = & \lambda\sqrt{\frac{t}{2\pi}}\int_{(2v+\epsilon)\sqrt{t}}^{\infty}e^{-(z-2v\sqrt{t})^{2}/2}dz\\
 & = & \lambda\sqrt{t}\bar{\Phi}(\epsilon\sqrt{t})\leq\lambda\sqrt{\frac{t}{2\pi}}e^{-\frac{\epsilon^{2}t}{2}}\end{eqnarray*}
 Thus one can see that for $\epsilon\sqrt{t}>1$ \begin{equation}
P(J_{\epsilon}(t))\leq\lambda\sqrt{\frac{2t}{\pi}}e^{-\frac{\epsilon^{2}t}{2}}\label{eq:est3}\end{equation}

As the event $A(t)$ belongs to the union of the events $C_{\epsilon}(t)$
and $J_{\epsilon}(t)$, then using the estimates (\ref{eq:est1})
and (\ref{eq:est3}), we get that for any sufficiently small $\epsilon>0$
\begin{equation}
P\left(N(t)>\frac{(\lambda+\epsilon)(v+\epsilon)t}{1-(\lambda+\epsilon)\delta}\right)=O\left(e^{-\alpha_{1}(\epsilon)t}\right)\label{eq:exp1}\end{equation}
 for some $\alpha_{1}(\epsilon)>0.$

Consider now the event $A'(t)=\{N(t,\omega)<c'_{\epsilon}t\},$ where
$c'_{\epsilon}=(1-(\lambda-\epsilon)\delta)^{-1}(\lambda-\epsilon)(v-\epsilon).$
Define the event $C'_{\epsilon}(t)$, that at time $0$ inside the
interval $(0,c'_{\epsilon}t\delta+(v-\epsilon)t)$ there were not
more than $(\lambda-\epsilon)(\delta c'_{\epsilon}t+(v-\epsilon)t)=c'_{\epsilon}t$
particles and let $\bar{C}'_{\epsilon}(t)$ be the complement to $C'_{\epsilon}(t).$
We can write the event $A'(t)$ as the union

\[
A'(t)=B_{\epsilon}'(t)\cup D_{\epsilon}'(t)\]
 where $B_{\epsilon}'(t)=A'(t)\cap C'_{\epsilon}(t)$ and $D_{\epsilon}'(t)=A'(t)\cap\bar{C}'_{\epsilon}(t).$
The event $B'_{\epsilon}(t)$ has exponentially small probability,
as it belongs to the event $C'_{\epsilon}(t),$ the probability of
which can be estimated via Poisson distribution:\begin{equation}
P(C'_{\epsilon}(t))\leq C'_{0}e^{-h'(\epsilon)t}\label{eq:est4}\end{equation}
 where the constant $C'_{0}$ depends on $\lambda$ and $\epsilon,$
$h'(\epsilon)>0.$

The event $D'_{\epsilon}(t)$ belongs to $J'_{\epsilon}(t),$ that
at least one particle of those initially in the interval $(0,c'_{\epsilon}t\delta+(v-\epsilon)t)$
was not absorbed by the boundary during time $t$. The probability
of event $J'_{\epsilon}(t)$ can be estimated from above by the expectation
of the number of particles, that were initially inside the interval
$(0,c'_{\epsilon}t\delta+(v-\epsilon)t),$ and deviated more than
on $\epsilon t$ its mean assignment point (if at time zero the particle
is at some point $x$, then its assignment point is $x-vt$) at time
$t$: \[
P(J'_{\epsilon}(t))\leq\lambda\int_{0}^{c'_{\epsilon}t\delta+(v-\epsilon)t}p'(t,x)dx\]
 where $p'(t,x)$ is the probability that a particle started from
$x$, $x<c'_{\epsilon}t\delta+(v-\epsilon)t,$ at time $t$ will be
outside the interval $(x-vt-\epsilon t,x-vt+\epsilon t),$ that is\[
p'(t,x)=1-\frac{1}{\sqrt{2\pi t}}\int_{x-(v+\epsilon)t}^{x-(v-\epsilon)t}e^{-\frac{(y-x+vt)^{2}}{2t}}dy\]
 After substituting $z=\frac{y-x+vt}{\sqrt{t}},$ we find\[
p'(t,x)=1-\frac{1}{\sqrt{2\pi}}\int_{-\epsilon\sqrt{t}}^{\epsilon\sqrt{t}}e^{-\frac{z^{2}}{2}}dy=2\bar{\Phi}(\epsilon\sqrt{t})\]
 Thus \[
P(J'_{\epsilon}(t))\leq C_{\epsilon}te^{-\frac{\epsilon^{2}t}{2}}\]
 Taking into account (\ref{eq:est4}), we get that for any sufficiently
small $\epsilon>0$\begin{equation}
P\left(N(t)<\frac{(\lambda-\epsilon)(v-\epsilon)t}{1-(\lambda-\epsilon)\delta}\right)=O\left(e^{-\alpha_{2}(\epsilon)t}\right)\label{eq:exp2}\end{equation}
 for some $\alpha_{2}(\epsilon)>0,$ together with (\ref{eq:exp1})
this proves the first assertion of the theorem.

\subsection{Zero drift }

Here we prove the second assertion of the Theorem.

We start with {}``trivial'' estimate from below, we mean by this
that we estimate just the number of particle collisions with the fixed
boundary. Consider an auxiliary model where the boundary does not
change its position while hit by particles but the particles disappear
at the first collision with the boundary. Let $L(t)$ be the number
of collisions with the boundary for this model during time $t$. Evidently,
we have $L(t)\leq N(t)$ a. s. Denote $K_{0}(t)=K_{0}(t,d)$ the number
of particles, which did not reach $0$ during time $t$ and at time
$t$ were to the left of the point $d\sqrt{t},$ where $d>0.$ Put
$K(t)=N(t)-L(t).$ If the condition $\delta L(t)\geq d\sqrt{t}$ holds
then $K(t)\geq K_{0}(t),$ as in the basic model all particles, situated
at time $t$ inside $[0,\delta L(t)],$ would have been absorbed by
the moving boundary. Thus one can state that\begin{equation}
MK(t)\geq MK_{0}(t)I(L(t)\geq d\sqrt{t})\label{eq:est_below}\end{equation}
 where $I(L(t)\geq d\sqrt{t})$ is the indicator of the corresponding
event.

Let us find the joint distribution of the variables $L(t)$ and $K_{0}(t).$
We begin with the means. Let $p(t,x)$ be the probability that a particle
situated initially at $x,$ during time $t$ will reach the boundary:\[
p(t,x)=P(\min_{s\leq t}(x+w(s))\leq0)\]
 As\begin{eqnarray*}
P(\min_{s\leq t}(x+w(s))\leq0) & = & P(\min_{s\leq t}w(s)\leq-x)\\
 & = & P(\max_{s\leq t}w(s)\geq x)\end{eqnarray*}
 then by formula (\ref{eq:distr_max}) with $v=0$ we have \[
p(t,x)=1-\Phi\left(\frac{x}{\sqrt{t}}\right)+\Phi\left(-\frac{x}{\sqrt{t}}\right)=2\bar{\Phi}\left(\frac{x}{\sqrt{t}}\right)\]
 Let us show that for zero drift ($v=0$)\[
ML(t)=\lambda\sqrt{\frac{2t}{\pi}}\]
 In fact, we have\begin{eqnarray*}
ML(t) & = & \lambda\int_{0}^{\infty}p(t,x)dx=2\lambda\int_{0}^{\infty}\bar{\Phi}\left(\frac{x}{\sqrt{t}}\right)dx\\
 & = & 2\lambda\sqrt{t}\int_{0}^{\infty}\bar{\Phi}\left(x\right)dx=\lambda\sqrt{\frac{2t}{\pi}}\end{eqnarray*}
 where we integrated the latter integral by parts\[
\int_{0}^{\infty}\bar{\Phi}\left(x\right)dx=\int_{0}^{\infty}xd\Phi(x)=\frac{1}{\sqrt{2\pi}}\int_{0}^{\infty}xe^{-x^{2}/2}dx=\frac{1}{\sqrt{2\pi}}\]

Let us find $MK_{0}(t).$ Denote $\sigma(t,x)$ the probability that
a particle, situated initially at $x,$ during time $t$ will not
reach $0$ and moreover at time $t$ it will be to the left of $d\sqrt{t}.$
Otherwise speaking\[
\sigma(t,x)=P(x+\min_{[0,t]}w_{s}>0,x+w_{t}<d\sqrt{t})=P(\max_{[0,t]}w_{s}<x,w_{t}>x-d\sqrt{t})\]
 The probability $\sigma(t,x)$ can be written as well in the following
way\[
\sigma(t,x)=P(w_{t}>x-d\sqrt{t})-2P(w_{t}>x)+P(\max_{[0,t]}w_{s}>x,w_{t}<x-d\sqrt{t})\]
 where due to \cite{Ito}, (see p. 45), we have\begin{eqnarray*}
P(\max_{[0,t]}w_{s}>x,w_{t}<x-d\sqrt{t}) & = & \sqrt{\frac{2}{\pi t^{3}}}\int_{-\infty}^{x-d\sqrt{t}}da\int_{x}^{\infty}(2b-a)e^{-(2b-a)^{2}/2t}db\\
 & = & \frac{1}{\sqrt{2\pi t}}\int_{-\infty}^{x-d\sqrt{t}}e^{-(2x-u)^{2}/2t}du\\
 & = & \bar{\Phi}\left(\frac{x+d\sqrt{t}}{\sqrt{t}}\right)\end{eqnarray*}
 Using this formula we find\[
\sigma(t,x)=\bar{\Phi}\left(\frac{x-d\sqrt{t}}{\sqrt{t}}\right)-2\bar{\Phi}\left(\frac{x}{\sqrt{t}}\right)+\bar{\Phi}\left(\frac{x+d\sqrt{t}}{\sqrt{t}}\right)\]
 and then\[
MK_{0}(t)=\int_{0}^{\infty}\sigma(t,x)dx=\sqrt{t}\left(\int_{-d}^{\infty}\bar{\Phi}(x)dx-2\int_{0}^{\infty}\bar{\Phi}(x)dx+\int_{d}^{\infty}\bar{\Phi}(x)dx\right)=\sigma\sqrt{t}\]
 where the constant \[
\sigma=\int_{-d}^{\infty}\bar{\Phi}(x)dx-2\int_{0}^{\infty}\bar{\Phi}(x)dx+\int_{d}^{\infty}\bar{\Phi}(x)dx=\int_{-d}^{0}\bar{\Phi}(x)dx-\int_{0}^{d}\bar{\Phi}(x)dx\]

\begin{lemma}\label{l_t}Random variables $L(t)$ and $K_{0}(t)$
are independent. $L(t)$ has Poisson distribution with parameter $\lambda\sqrt{\frac{2t}{\pi}},$
and $K_{0}(t)$ has Poisson distribution with parameter $\sigma\sqrt{t}.$

\end{lemma}

Proof. Let $\pi_{n}(a)$ be Poisson-distribution with parameter $\lambda a$,
$L(a,t)$ be the number of particles, which were inside $(0,a)$ at
time $t=0$ and during $t$ were absorbed by the boundary, $K_{0}(a,t)$
be the number of particles which were inside $(0,a)$ at time $t=0$,
during time $t$ did not reach $0$ and at time $t$ were to the left
of the point $d\sqrt{t}.$ Put\[
q(t,x)=1-p(t,x)-\sigma(t,x)\]
 Find the joint distribution of the random variables $L(a,t)$ and
$K_{0}(a,t):$\[
P(L(a,t)=k,K_{0}(a,t)=m)=\sum_{n=k+m}^{\infty}\pi_{n}(a)a^{-n}\sum\int_{0}^{a}\dots\int_{0}^{a}\prod_{l=1}^{k}p(t,x_{i_{l}})dx_{i_{l}}\prod_{l=1}^{m}\sigma(t,x_{j_{l}})dx_{j_{l}}\prod_{l=1}^{n-k-m}q(t,x_{r_{l}})dx_{r_{l}}\]
 where the inner sum is over all ordered arrays $i_{1}<\dots<i_{k}$
and $j_{1}<\dots<j_{m}$, of length $k$ and $m$ correspondingly,
such that $\{i_{1}<\dots<i_{k}\}\cap\{j_{1}<\dots<j_{m}\}=\emptyset.$
Here\[
\{r_{1}<\dots<r_{n-k-m}\}=\{1,...,n\}\setminus\{i_{1}<\dots<i_{k}\}\cup\{j_{1}<\dots<j_{m}\}\]
 As the inner sum consists of \[
\frac{n!}{k!m!(n-k-m)!}\]
 equal terms $\left(\hat{p}(t,a)\right)^{k}\left(\hat{\sigma}(t,a)\right)^{m}\left(\hat{q}(t,a)\right)^{n-k-m},$
where $\hat{p}(t,a)=\int_{0}^{a}p(t,x)dx,$ $\hat{\sigma}(t,a)=\int_{0}^{a}\sigma(t,x)dx$
, $\hat{q}(t,a)=\int_{0}^{a}q(t,x)dx$, then\begin{eqnarray*}
P(L(a,t)=k,K_{0}(a,t)=m) & = & \sum_{n=k+m}^{\infty}\pi_{n}(a)a^{-n}\frac{n!}{k!m!(n-k-m)!}\left(\hat{p}(t,a)\right)^{k}\left(\hat{\sigma}(t,a)\right)^{m}\left(\hat{q}(t,a)\right)^{n-k-m}\\
 & = & \frac{\left(\lambda\hat{p}(t,a)\right)^{k}}{k!}\frac{\left(\lambda\hat{\sigma}(t,a)\right)^{m}}{m!}e^{-\lambda a}\sum_{n=k+m}^{\infty}\frac{\left(\lambda\hat{q}(t,a)\right)^{n-k-m}}{(n-k-m)!}\\
 & = & \frac{\left(\lambda\hat{p}(t,a)\right)^{k}e^{-\lambda\hat{p}(t,a)}}{k!}\frac{\left(\lambda\hat{\sigma}(t,a)\right)^{m}e^{-\lambda\hat{\sigma}(t,a)}}{m!}\end{eqnarray*}
 To get joint distribution of $L(t)$ and $K_{0}(t),$ let us tend
$a$ to infinity. In the limit we will get direct product of Poisson
distributions with parameters $\lambda\sqrt{\frac{2t}{\pi}}$ and
$\sigma\sqrt{t},$ as $\hat{p}(t,a)\to\sqrt{\frac{2t}{\pi}}$ and
$\hat{\sigma}(t,a)\to\sigma\sqrt{t}$ for $a\to\infty.$ The lemma
is proved.

Let us prove now that the constant in the lower bound is in fact greater
than the trivial one. Choose $d>0$ so that $d<\delta ML(t)/\sqrt{t}.$
By (\ref{eq:est_below}) and independence of random variables $L(t)$
and $K_{0}(t),$ we have \[
MK(t)\geq MK_{0}(t)P(L(t)\geq d\sqrt{t})\]
 and\[
MK_{0}(t)P(L(t)\geq d\sqrt{t})=\sigma\sqrt{t}-MK_{0}(t)P(L(t)<d\sqrt{t})\]
 As, for a given choice of $d$, the probability of the event $L(t)<d\sqrt{t}$
is exponentially small\[
P(L(t)<d\sqrt{t})\leq C_{4}e^{-\phi\sqrt{t}}\]
 then one can assert that $MK_{0}(t)P(L(t)\geq d\sqrt{t})\geq\sigma_{0}\sqrt{t}$
for some constant $\sigma_{0}>0,$ where $\sigma_{0}<\sigma.$

Prove now the bound from above in the second statement of the Theorem.
Our goal is to estimate the probability $P(N(t)>x\sqrt{t})$ from
above. Take arbitrary $\epsilon>0,$ so that $(\lambda+\epsilon)\delta<1.$
Subdivide the event $A(t,N)=\{N(t,\omega)=N\}$ into the union of
events\[
A(t,N)=\cup_{n\leq N}B_{0}(N,n,t)\]
 where $B_{0}(N,n,t)$ is the event that the boundary absorbed exactly
$n$ particles of the particles, situated initially inside $(0,N\delta)$,
and exactly $m=N-n$ particles of the particles, situated initially
inside $(N\delta,\infty)$. Denote $B_{0,\epsilon}(N,t)$ the union
of events $B_{0}(N,n,t)$ with $n>(\lambda+\epsilon)\delta N$. The
event $B_{0,\epsilon}(N,t)$ has exponentially small probability,
as this event belongs to the event $C_{0,\epsilon}(N)$, that at time
$0$ inside $(0,N\delta)$ there was no more than $(\lambda+\epsilon)\delta N$
particles. The probability of the latter event can be estimated by
the Poisson distribution:\begin{equation}
P(C_{0,\epsilon}(N))\leq C_{0}e^{-h(\epsilon)N}\label{eq:e1}\end{equation}
 where $h(\epsilon)>0$ and $C_{0}$ is some constant depending only
on $\lambda,\epsilon.$

Denote $D_{0,\epsilon}(N,t)$ the union of events $B_{0}(N,n,t)$
with $n<(\lambda+\epsilon)\delta N$. From the considerations above
it follows that the event $D_{0,\epsilon}(N,t)$ is equivalent to
the event \[
m=N-n>N(1-(\lambda+\epsilon)\delta)\]
 At the same time the event $m>N(1-(\lambda+\epsilon)\delta)$ belongs
to the event that among the particles which were at time $0$ to the
right of the point $N\delta,$ more than $N(1-(\lambda+\epsilon)\delta)$
particles during time $t$ were at least once to the left of the point
$N\delta$. Denote the latter event by $J_{0,\epsilon}(N,t)$, and
denote $\kappa(t)$ the number of particles which were at time $0$
to the right of the point $N\delta$ but having visited the left side
of $N\delta$ at least once during time $t$ (that is $P(J_{0,\epsilon}(N,t))=P(\kappa(t)>N(1-(\lambda+\epsilon)\delta)$).
The distribution of this random variable does not depend on $N$ and
coincides with the distribution of the random variable $L(t)$, defined
above. By lemma \ref{l_t} and the properties of Poisson distribution,
we have the following estimate \begin{equation}
P(\kappa(t)>y\sqrt{t})\leq C_{1}e^{-\beta\sqrt{t}}\label{eq:e2}\end{equation}
 valid for some $\beta>0$, $y\geq y_{0}>\lambda\sqrt{\frac{2}{\pi}}$
and constant $C_{1},$ depending on $y_{0}.$

As the event $A(t,N)$ belongs to the union of the events $C_{0,\epsilon}(N)$
and $J_{0,\epsilon}(N,t),$ then \[
P(N(t)>x\sqrt{t})\leq\sum_{N>x\sqrt{t}}P(C_{0,\epsilon}(N))+\sum_{N>x\sqrt{t}}P(J_{0,\epsilon}(N,t))\]
 By (\ref{eq:e1}), for the first sum the following estimate holds
\[
\sum_{N>x\sqrt{t}}P(C_{0,\epsilon}(N))=O\left(e^{-h(\epsilon)x\sqrt{t}}\right)\]
 The second sum, by (\ref{eq:e2}), has the following estimate\[
\sum_{N>x\sqrt{t}}P(J_{0,\epsilon}(N,t))=\sum_{N>x\sqrt{t}}P(\kappa(t)>(1-(\lambda+\epsilon)\delta)N)\leq C_{2}e^{-\beta(\epsilon)\sqrt{t}}\]
 for some constant $\beta(\epsilon)>0$, $x>x_{0}=y_{0}(1-(\lambda+\epsilon)\delta)^{-1}$
and some constant $C_{2},$ depending on $x_{0}.$ It follows, that
for sufficiently large $x$ \[
P(N(t)>x\sqrt{t})\leq C_{3}e^{-\gamma x\sqrt{t}}\]
 Hence, for $t$ sufficiently large\[
MN(t)\leq C\sqrt{t}\]

The Theorem is proved.

\paragraph{Proof of the second part of Lemma 1}

Consider first the case of zero drift. Let us find the distribution
of particle configuration at the time of first hitting the boundary.
Put\[
\sigma(t,x,y)dy=P(x+w_{t}\in dy,x+\min_{s\leq t}w_{t}>0),x>0,y>0\]
 Thus, $\sigma(t,x,y)dy$ is the probability that the particle which
were at time $t=0$ at $x>0,$ during time $t$ will not reach $0$
and will be in the interval $(y,y+dy)$ at time $t$. As the distribution
of $w_{t}$ coincides with the distribution of $-w_{t},$ we have\[
\sigma(t,x,y)dy=P(x-w_{t}\in dy,x-\max_{s\leq t}w_{t}>0)=P(w_{t}\in d(x-y),\max_{s\leq t}w_{t}<x)\]
 To calculate $\sigma(t,x,y)$ we will use the following formula for
the joint distribution (see, for example \cite{Ito})\[
P(w_{t}\in da,\max_{s\leq t}w_{s}\in db)=\sqrt{\frac{2}{\pi t^{3}}}(2b-a)e^{-(2b-a)^{2}/2t}dadb,0\leq b,b\geq a\]
 The change of variables $u=(2b-x+y)^{2}/2t$ gives \begin{eqnarray*}
\sigma(t,x,y) & = & \sqrt{\frac{2}{\pi t^{3}}}\int_{\max(x-y,0)}^{x}(2b-x+y)e^{-(2b-x+y)^{2}/2t}db\\
 & = & \sqrt{\frac{1}{2\pi t}}\int_{(x-y)^{2}/2t}^{(x+y)^{2}/2t}e^{-u}du\\
 & = & \sqrt{\frac{1}{2\pi t}}\left(e^{-(x-y)^{2}/2t}-e^{-(x+y)^{2}/2t}\right)\end{eqnarray*}
 Let $\tau$ be the moment of first hitting the boundary by one of
the particles.

\begin{lemma}The probability density of the random variable $\tau$is\[
P(\tau\in dt)=\frac{\lambda}{\sqrt{2\pi t}}e^{-\lambda\sqrt{2t/\pi}}dt\]
 Under the condition $\tau=t$ the particle configuration at time
$t$ is an inhomogeneous Poisson point field with the rate\[
\lambda\psi(t,y)=\lambda\int_{0}^{\infty}\sigma(t,x,y)dx=\lambda(\Phi(y/\sqrt{t})-\Phi(-y/\sqrt{t}))\]
 where $\Phi(y)=\frac{1}{\sqrt{2\pi}}\int_{-\infty}^{y}e^{-z^{2}/2}dz.$

\end{lemma}

Proof. Let $B_{1},...,B_{l}$ be an array consisting of $l$ pairwise
non intersecting intervals on the positive axis. Let us find the joint
distribution $P(\tau\in dt,\eta(t,B_{i})=k_{i},i=1,...l),$ where
$\eta(t,B_{i})$ is the number of particles in the interval $B_{i}$
at time $t.$ Let $\pi_{n}(a)$ be the Poisson distribution with parameter
$\lambda a$, $\eta(a,t,B_{i})$ be the number of particles, which
were initially inside the interval $(0,a)$ and at time $t$ were
in $B_{i}.$ Let $\tau(a)$ be the first time moment, when one of
the particles, which were in $(0,a)$ at time $t=0,$ hits the boundary.
Denote $g(t,x)$ the probability density of the time $\beta_{x}$
of first hitting the point $0$ by the particle under the condition
that at time $t=0$ this particle were at $x.$ It is well known that
(see for example \cite{Ito}), that\[
g(t,x)=\frac{x}{\sqrt{2\pi t^{3}}}e^{-x^{2}/2t}\]
 Let us find the joint distribution of random variables $\eta(a,t,B_{i}),i=1,...,l$
and $\tau(a):$\begin{eqnarray*}
P(\tau(a)\in dt,\eta(a,t,B_{i})=k_{i},i=1,...l) & = & \sum_{n=1+k_{1}+...+k_{l}}^{\infty}\pi_{n}(a)a^{-n}\sum\int_{0}^{a}\dots\int_{0}^{a}g(t,x_{m})dtdx_{m}\\
 &  & \prod_{j=1}^{l}\prod_{s=1}^{k_{j}}dx_{i_{j,s}}\int_{y\in B_{j}}\sigma(t,x_{i_{j,s}},y)dy\\
 &  & \prod_{s=1}^{n-1-k_{1}-...-k_{l}}dx_{r_{s}}\int_{y\in R_{+}\setminus\cup_{j=1}^{l}B_{j}}\sigma(t,x_{r_{s}},y)dy\end{eqnarray*}
 where the inner sum is over all pairwise non intersecting ordered
arrays $m,$ $i_{1,1}<...<i_{1,k_{1}},$ $i_{2,1}<...<i_{2.k_{2}},$...,
$i_{l,1}<...<i_{l.k_{l}}$ of the lengths $1,k_{1},k_{2},...k_{l}$
correspondingly. Here\[
\{r_{1}<\dots<r_{n-1-k_{1}-...-k_{l}}\}=\{1,...,n\}\setminus(\{m\}\cup\{i_{1,1}<...<i_{1,k_{1}}\}\cup...\cup\{i_{l,1}<...<i_{l.k_{l}}\})\]
 As the inner sum consists of\[
\frac{n!}{k_{1}!...k_{l}!(n-1-k_{1}-...-k_{l})!}\]
 equal terms $\hat{g}(a,t)dt\prod_{j=1}^{l}\left(\hat{\sigma}(a,t,B_{j})\right)^{k_{j}}\left(\hat{\sigma}(a,t,R_{+}\setminus\cup_{j=1}^{l}B_{j})\right)^{n-1-k_{1}-...-k_{l}},$
where\begin{eqnarray*}
\hat{g}(t,a) & = & \int_{0}^{a}g(t,x)dx\\
\hat{\sigma}(a,t,B_{j}) & = & \int_{0}^{a}dx\int_{y\in B_{j}}\sigma(t,x,y)dy\end{eqnarray*}
 then\begin{eqnarray*}
P(\tau(a)\in dt,\eta(a,t,B_{i})=k_{i},i=1,...l) & = & \sum_{n=1+k_{1}+...+k_{l}}^{\infty}\pi_{n}(a)a^{-n}\frac{n!}{k_{1}!...k_{l}!(n-1-k_{1}-...-k_{l})!}\\
 &  & \hat{g}(a,t)dt\prod_{j=1}^{l}\left(\hat{\sigma}(a,t,B_{j})\right)^{k_{j}}\left(\hat{\sigma}(a,t,R_{+}\setminus\cup_{j=1}^{l}B_{j})\right)^{n-1-k_{1}-...-k_{l}}\\
 & = & \lambda\hat{g}(a,t)dt\prod_{j=1}^{l}\frac{\left(\lambda\hat{\sigma}(a,t,B_{j})\right)^{k_{j}}}{k_{j}!}e^{-\lambda a}\times\\
\\ &  & \times\sum_{n=1+k_{1}+...+k_{l}}^{\infty}\frac{\left(\lambda\hat{\sigma}(a,t,R_{+}\setminus\cup_{j=1}^{l}B_{j})\right)^{n-1-k_{1}-...-k_{l}}}{(n-1-k_{1}-...-k_{l})!}\\
 & = & \lambda\hat{g}(a,t)dt\prod_{j=1}^{l}\frac{\left(\lambda\hat{\sigma}(a,t,B_{j})\right)^{k_{j}}}{k_{j}!}e^{-\lambda a+\lambda\hat{\sigma}(a,t,R_{+}\setminus\cup_{j=1}^{l}B_{j})}\end{eqnarray*}
 Taking into account that\[
\sum_{j=1}^{l}\hat{\sigma}(a,t,B_{j})+\hat{\sigma}(a,t,R_{+}\setminus\cup_{j=1}^{l}B_{j})=\int_{0}^{a}P(\beta_{x}>t)dx=a-\int_{0}^{a}P(\beta_{x}<t)dx\]
 we have\[
-a+\hat{\sigma}(a,t,R_{+}\setminus\cup_{j=1}^{l}B_{j})=-\sum_{j=1}^{l}\hat{\sigma}(a,t,B_{j})-\int_{0}^{a}P(\beta_{x}<t)dx\]
 Thus\[
P(\tau(a)\in dt,\eta(a,t,B_{i})=k_{i},i=1,...l)=\lambda\hat{g}(a,t)e^{-\int_{0}^{a}P(\beta_{x}<t)dx}dt\prod_{j=1}^{l}\frac{\left(\lambda\hat{\sigma}(a,t,B_{j})\right)^{k_{j}}}{k_{j}!}e^{-\lambda\hat{\sigma}(a,t,B_{j})}\]
 To obtain the joint distribution of $\tau$ and $\eta(t,B_{i}),$
we tend $a$ to infinity. Then\begin{eqnarray*}
\hat{g}(a,t) & \to & \int_{0}^{\infty}g(t,x)dx=\sqrt{\frac{1}{2\pi t}}\\
\int_{0}^{a}P(\beta_{x}<t)dx & \to & \int_{0}^{\infty}P(\beta_{x}<t)dx=\sqrt{\frac{2t}{\pi}}\\
\hat{\sigma}(a,t,B_{j}) & \to & \hat{\sigma}(t,B_{j})=\int_{y\in B_{j}}dy\int_{0}^{\infty}\sigma(t,x,y)dx\\
 & = & \int_{y\in B_{j}}dy\psi(t,y)\end{eqnarray*}
 where $\psi(t,y)=\Phi(y/\sqrt{t})-\Phi(-y/\sqrt{t}).$ Hence\[
P(\tau\in dt,\eta(t,B_{i})=k_{i},i=1,...l)=\frac{\lambda}{\sqrt{2\pi t}}e^{-\lambda\sqrt{\frac{2t}{\pi}}}dt\prod_{j=1}^{l}\frac{\left(\lambda\hat{\sigma}(t,B_{j})\right)^{k_{j}}}{k_{j}!}e^{-\lambda\hat{\sigma}(t,B_{j})}\]
 The lemma is proved.

Let us prove that the random variable $k{}_{1}$ can take infinite
values with positive probability. Assume that $\tau=t.$ Denote $m_{1}$
the number of particles in the interval $[0,\delta]$ at time $\tau=t$,
$m_{2}$ the number of particles at time $\tau=t$ in the interval
$D_{1}=[\delta,(m_{1}+1)\delta]$, if $m_{1}>0$, and denote $m_{3}$
the number of particles at time $\tau=t$ in the interval $D_{2}=[(m_{1}+1)\delta,(m_{1}+1)\delta+m_{2}\delta],$
if $m_{1}>0$, $m_{2}>0$, etc. The sequence of pairs $\mu_{i}=(m_{i-1},m_{i})$
is a discrete time homogeneous Markov chain. Transition probabilities
are\[
p_{\mu_{i}\mu_{i+1}}=\frac{(\lambda\hat{\sigma}(t,D_{i}))^{m_{i+1}}}{m_{i+1}!}e^{-\lambda\hat{\sigma}(t,D_{i})}\]
 The states $(m,0)$ are absorbing. Random variable $k{}_{1}$ is
defined by the random stopping time of the Markov chain $\mu_{i}$.
The boundary stops, when the chain hits some of the absorbing states.
We introduce Lyapounov function $f(\mu_{i})=m_{i}.$ Its one step
increment equals\begin{eqnarray*}
M(f(\mu_{i+1})-f(\mu_{i})|\mu_{i}=(m_{i-1},m_{i})) & = & \lambda\hat{\sigma}(t,D_{i})-m_{i}\\
 & = & \lambda\int_{(m_{i-1}+1)\delta}^{(m_{i-1}+1)\delta+m_{i}\delta}\psi(t,y)dy-m_{i}\end{eqnarray*}
 As $\psi(t,y)\to1$ for $y\to\infty$ and $\lambda\delta>1,$ then
there exists such $\epsilon>0,$ that for sufficiently large $m_{i-1}$
\[
M(f(\mu_{i+1})-f(\mu_{i})|\mu_{i}=(m_{i-1},m_{i}))\geq\epsilon m_{i}\]
 This is the transience of the Markov chain $\mu_{i}.$ Thus, $P(k_{1}=\infty)>0.$

The case of nonzero drift is considered similarly. Moreover, as we
consider the lower bound, it is intuitively clear that the drift can
only increase the probability that the boundary explodes to infinity.


\begin{thebibliography}{3}
\bibitem{Shir} A. N. Shiryaev. Foundations of stochastic financial
mathematics. V. 2, Moscow 1998.

\bibitem{Ito} K. Ito, G. McKean. Diffusion processes and their trajectories.
1968.

\bibitem{Dobr}R. L. Dobrushin. About Poisson law of particle distribution
in space. Ukrainian Math. Journal, 1956, v. 8, No. 2, 130-134. 
\end{thebibliography}
\end{document}